\documentclass[11pt]{article}
\usepackage{amsmath,amsthm,amsfonts,amssymb,graphicx}

\title{On the residual finiteness of outer automorphisms of relatively hyperbolic groups}
\author{V. Metaftsis and M. Sykiotis}

\makeatletter
\date{}
\newcommand{\lar}{\rightarrow}
\newcommand{\dsp}{\displaystyle}
\newtheorem{thm}{Theorem}[section]
\newtheorem{lem}[thm]{Lemma}
\newtheorem{cor}[thm]{Corollary}

\theoremstyle{definition}

\theoremstyle{remark}

\def\classification{\@ifnextchar [{\@xfootnotetext}%
{\begingroup\let\protect\noexpand\xdef\@thefnmark{}
\endgroup\@footnotetext}}
\def\keywords{\@ifnextchar [{\@xfootnotetext}%
{\begingroup\let\protect\noexpand\xdef\@thefnmark{}
\endgroup\@footnotetext}}
\newcommand{\aut}{{\rm Aut}}
\newcommand{\out}{{\rm Out}}
\newcommand{\inn}{{\rm Inn}}
\newcommand{\conj}{{\rm Conj}}
\renewcommand{\mod}{\rm Mod}
\newcommand{\R}{\mathbb R}

\newcommand{\N}{\mathbb N}
\renewcommand{\H}{\mathbb H}
\renewcommand{\t}{\tau}
\renewcommand{\d}{\delta}
\renewcommand{\l}{\lambda}
\newcommand{\w}{\omega}
\begin{document}
\classification {2000 {\sl Mathematics Subject Classification.}
20F67, 20E26, 20E36, 20F28, 20E07, 57M60}
\maketitle

\begin{abstract} We show that every virtually torsion-free subgroup of the outer automorphism group
of a conjugacy separable hyperbolic group is residually finite. As a result, we are able to prove 
that the group of outer automorphisms of every finitely generated Fuchsian group and of every
free-by-finite group is residually 
finite. In an addendum, we also generalise the main result for relatively hyperbolic groups.
\end{abstract}

\section{Introduction}
The {\it mapping class group} of a compact connected, possibly with boundary, surface  
$S$, $\mod_S$, is defined as 
the group of isotopy
classes of homeomorphisms $S\lar S$. The algebraic investigation
of that group goes as far back as Dehn and Nielsen, who first showed that the mapping class
group of an orientable surface is embedded into the group of outer automorphisms of the
fundamental group of $S$, $\out(\pi_1(S))$. This embedding was later extended to a wider 
class of 2-orbifolds. Maclachlan and Harvey \cite{mh} proved that $\mod_S\le \out(\pi_1(S))$ 
for a surface $S$ obtained from a compact surface of genus $k$ by deleting $s$ points,
$t$ discs and $r$ marked points. Recently, Fujiwara \cite{fujiwara} extended the above 
embedding to the case of hyperbolic 2-orbifolds (orientable or non-orientable) with finite volume. 

Grossman \cite{Gr} was the first to show that if $S_k$ is a compact
orientable surface of genus $k$, then $\mod_{S_k}$
is residually finite. Recently, Allenby, Kim and Tang \cite{akt} extended the result
of Grossman for non-orientable closed surfaces. Ivanov \cite{ivanov} gave a geometric proof
of the result of Grossman which seems that can easily be extended to the non-orientable 
case. The proofs of both results in \cite{Gr} and \cite{akt} are using
combinatorial group theory arguments. One of the key ingredients of both results
is the fact that the group of conjugating automorphisms of the fundamental group of the surface
coincides with the group of its inner automorphisms. If $G$ is a group, then an automorphism 
$f$ of $G$ is called a {\it conjugating automorphism} if  $f(g)$ is a conjugate of $g$ for 
every $g\in G$. The conjugating automorphisms of a group $G$ form a  subgroup of the group 
of  automorphisms of $G$, $\aut(G)$, which
we denote by $\conj(G)$. Moreover, $\conj(G)$ is normal in $\aut(G)$ and 
contains the subgroup of inner automorphisms of $G$, $\inn(G)$.

In the present note we investigate the residual finiteness of $\out(G)$ for hyperbolic
groups $G$. Our study is, as well, based on the investigation of the group of conjugating 
automorphisms of $G$.
In fact, using the powerful geometric ideas developed by Paulin \cite{Pau} and 
Bestvina \cite{bestvina} in the context of ultralimits of metric spaces, 
we show that in a hyperbolic group $G$, the quotient group $\conj(G)/\inn(G)$ is
always finite. The main theorem shows that every virtually torsion-free subgroup of the
outer automorphism group of a conjugacy separable, hyperbolic group is residually
finite.

As a result, we obtain the residual finiteness of the outer automorphism group 
of finitely generated Fuchsian groups and of free-by-finite groups. Consequently the 
mapping class groups of hyperbolic 2-orbifolds with finite volume are 
residually finite. Hence, we retrieve the results of Grossman and of
Allenby, Kim and Tang as special cases of our corollaries.

In an addendum at the end of the paper we show how to generalize the main
results for relatively hyperbolic groups.

\section{Main Results}

The next lemma is what really proved in ~\cite[Theorem 1]{Gr}. We
reproduce here the proof for the reader's convenience.

\begin{lem}\label{lem:Gr} Let $G$ be a finitely generated, conjugacy separable
group. Then the quotient group $\aut(G)/\conj(G)$ is residually
finite.
\end{lem}

\begin{proof}
Since $G$ is finitely generated, we can find a family
$(K_{i})_{i\in I}$ of characteristics subgroups of finite index in
$G$ such that for each normal subgroup $N$ of finite index in $G$
there is a subgroup $K_{i}$ in the above family, contained in $N$. 
Every automorphism $f$ of $G$ induces an
automorphism $f_{i}$ of $G/K_{i}$ which acts as a permutation on the
conjugacy classes of $G/K_{i}$. Therefore, for each $i\in I$ we have a
homomorphism $\pi_{i}:\aut(G)\lar S_{i}$, where $S_{i}$ denotes the
permutation group on the set of conjugacy classes of the finite group
$G/K_{i}$. Since each conjugating automorphism preserves conjugacy
classes, the group $\conj(G)$ is contained in the intersection
$\bigcap _{i\in I}$ Ker$(\pi_{i})$ of the kernels of the $\pi_{i}$'s.
Now let $f\in \bigcap _{i\in I}$ Ker$(\pi_{i})$. Then for every $g\in G$, 
the elements $f(g)$ and $g$ are conjugate in each quotient
$G/K_{i}$. In particular, the elements $f(g)$ and $g$ are
conjugate in each finite quotient of $G$ and thus they are
conjugate in $G$, by the conjugacy separability of $G$. This means
that $f$ is a conjugating automorphism of $G$ and hence
$\conj(G)=\bigcap_{i\in I}$ Ker$(\pi_{i})$.
\end{proof}

The following material on ultralimits of metric spaces can be
found in detail in the paper of Kapovich and Leeb \cite{kl}.

A $\textsl{filter}$ on the set of natural numbers $\mathbb{N}$ is
a non-empty set $\omega$ of subsets of $\mathbb{N}$ with the
following properties:

\begin{enumerate}
\item The empty set is not contained in $\omega$.
\item $\omega$ is closed under finite intersections.
\item If $A\in \omega$ and $A\subseteq B$, then $B\in \omega$.
\end{enumerate}

A filter $\omega$ is called $\textsl{ultrafilter}$ if it is
maximal. An ultrafilter is called \textsl{non-principal} if it
contains the complements of the finite subsets of $\mathbb{N}$.
Let $\omega$ be an ultrafilter on $\mathbb{N}$ and
$a:\mathbb{N}\lar [0,\infty)$ a sequence of non-negative real
numbers. Then there is a unique point $l$, which we denote by
$\textrm{lim}_{\omega}a_{i}$, in the one-point compactification
$[0,\infty]$ of $[0,\infty)$ such that for each neighborhood $U$
of $l$, the inverse image $a^{-1}(U)$ of $U$ under $a$ is
contained in $\omega$.

Let $(X_{i},d_{i},x_{i}^{0})_{i\in \mathbb{N}}$ be a sequence of
based metric spaces and let $\omega$ be an ultrafilter on
$\mathbb{N}$. On the subspace $Y$ of $\prod_{i\in
\mathbb{N}}X_{i}$ consisting of all sequences $(x_{i})$ for which
$\textrm{lim}_{\omega}d_{i}(x_{i},x_{i}^{0})<\infty$, we define a
pseudo-metric $d_{\omega}$ by
$d_{\omega}\big((x_{i}),(y_{i})\big)=\textrm{lim}_{\omega}\big(d_{i}(x_{i},y_{i})\big)$.
The \textsl{based ultralimit}
$(X_{\omega},d_{\omega},x_{\omega})$, where
$x_{\omega}=(x_{i}^{0})_{i\in \mathbb{N}}$, is the associated
metric space.

\bigskip

Let $X$ be a geodesic metric space and $\delta$ a non-negative real number.
The space $X$ is called {\it  $\delta$-hyperbolic\/} if for every triangle $\Delta\subset X$
with geodesic sides, each side is contained in the $\delta$-neighbourhood of the
union of the two other sides.

Let $G$ be a finitely generated group and let $X=X(G,S)$ be the Cayley graph of
$G$ with respect to a finite generating set $S$ for $G$, closed under inverses.
The group $G$ is called {\it $\delta$-hyperbolic\/} if $X$ is a $\delta$-hyperbolic
space with respect to the word metric.

In the next lemma we use various results concerning actions of groups on $\R$-trees.
For details, we refer the reader 
to the paper of Morgan and Shalen \cite{ms}.

The key point of Lemma \ref{lem:finite} is the following statement.
Assume that an action of a hyperbolic group $G$ on a real tree $Y$ is obtained 
as a limit of a sequence of actions of $G$ on its Cayley graph $X$, where each 
action in the sequence is the natural action of $G$ on $X$ twisted by an 
automorphism of $G$. If the conjugacy class of an element $g$ in $G$ is periodic 
under these automorphisms, then $g$ is elliptic when acting on the limit tree $Y$.

As it was pointed out by the referee, the above statement seems
to be well known to the experts (see \cite{bestvina2}). Nonetheless, the
authors failed to track down a reference for its proof.  So, a proof of it,
is included in the lemma for the reader's convenience and completeness. 

\begin{lem}\label{lem:finite} Let $G$ be a hyperbolic group. Then the group 
$\inn(G)$ of inner automorphisms of $G$ is of
finite index in $\conj(G)$.
\end{lem}

\begin{proof} Suppose to the contrary that $\inn(G)$ is of infinite
index in $\conj(G)$. Fix an infinite sequence
$f_{1},f_{2},\dots,f_{n},\dots$ of conjugating automorphisms of
$G$ representing pairwise distinct cosets of $\inn(G)$ in $\conj(G)$. We
apply the method of Bestvina and Paulin, using ultrafilters, to
construct an $\mathbb{R}$-tree on which $G$ acts by isometries.

Let $X=X(G,S)$ be the Cayley graph of $G$ with respect to a finite
generating set $S$ closed under inverses. Then $X$ with the associated
word metric $d$ is a $\delta$-hyperbolic metric space for some $\delta\ge 0$. 
Each $f_{i}$ gives an
action $\rho_{i}$ by isometries of $G$ on $X$ by
$\rho_i(g,x)=f_{i}(g)x$. The outer automorphism group of a
virtually cyclic group is finite. Thus, we may assume that $G$ is
non-elementary. In that case the action of $G$ on the boundary
$\partial X$ of $X$ is non-trivial and Lemma 2.1 in ~\cite{Pau}
applies to any action $\rho_{i}$, yielding a sequence $x_{i}^{0}$
of elements of $X$ such that
\begin{equation} \label{eq:min} \underset{g\in
S}{\textrm{max}}\,d\big(x_{i}^{0},f_{i}(g)x_{i}^{0}\big)\leq
\underset{g\in S}{\textrm{max}}\,d\big(x,f_{i}(g)x\big)
\end{equation}
for all $x$ in $X$. Let $\lambda_{i}=\underset{g\in
S}{\textrm{max}}\,d\big(x_{i}^{0},f_{i}(g)x_{i}^{0}\big)$. Then
\begin{equation} \label{eq:trineq}
d\big(x_{i}^{0},f_{i}(g)x_{i}^{0}\big)\leq \lambda_{i} \|g\|
\end{equation}
for all $g$ in $G$, where $\|g\|$ denotes the word-length of $g$
with respect to $S$. If the sequence $(\lambda_{i})$ contains a
bounded subsequence, then the argument in ~\cite[Case 1, p.
338]{Pau} shows that there are indices $i$ and $j$ with $i\neq j$
such that the automorphisms $f_{i}$ and $f_{j}$ differ by an inner automorphism
of $G$ and consequently they give rise to the same coset of $\inn(G)$, 
contradicting the choice of the $f_{i}$. It
follows that
$\underset{i\rightarrow\infty}{\textrm{lim}}\lambda_{i}=\infty$.
We consider the sequence $(X_{i},d_{i},x_{i}^{0})$ of based metric
spaces, where $X_{i}=X$ and $d_{i}=\frac{d}{\lambda_{i}}$. Note
that the space $(X_{i},d_{i},x_{i}^{0})$ is
$\frac{\delta}{\lambda_{i}}$-hyperbolic for all $i$.
For any non-principal
ultrafilter $\omega$ on $\mathbb{N}$, let
$(X_{\omega},d_{\omega},x_{\omega})$ be the corresponding based
ultralimit. The fact that the distance $d_{\omega}(x,y)$ of two
points $x=(x_{i})$ and $y=(y_{i})$ of $X_{\omega}$ is approximated
by the distances $d_{i}(x_{i},y_{i})$ for infinitely many $i$,
implies that $X_{\omega}$ is a $0$-hyperbolic space (i.e., an
$\mathbb{R}$-tree), since $\underset{i\rightarrow
\infty}{\textrm{lim}}\frac{\delta}{\lambda_{i}}=0$. The action of
$G$ on $X_{\omega}$ is given by
$g\cdot(x_{i})=\big(f_{i}(g)x_{i}\big)$. Inequality
~(\ref{eq:trineq}) ensures that the action is well-defined. We will
show that the action of $G$ on $X_{\omega}$ is trivial (i.e. there
is a global fixed point).

Let $g$ be an element of $G$ which acts as a hyperbolic isometry
on $X_{\omega}$, and let $\t_{\omega}(g)$ denote its translation
length. Fix $x=(x_{i})\in X_{\omega}$ such that $x$ lies on the axis of $g$. Then
$$\t_{\omega}(g)=d_{\omega}(gx,x)=\left(\textrm{lim}_{\omega}
d_i(f_i(g)x_i,x_i)=\right)
\dsp\frac{d_{\omega}(g^{n}x,x)}{n}$$ for all positive integers $n$.
In particular,
$d_{\omega}(gx,x)=\dsp\frac{d_{\omega}(g^{2}x,x)}{2}$ and thus
\begin{equation}\label{eq:2}
\textrm{lim}_{\omega}\Big(2d_{i}\big(f_{i}(g)x_{i},x_{i}\big)-d_{i}
\big(f_{i}(g)^{2}x_{i},x_{i}\big)\Big)
=0.
\end{equation}
For each $i$, we fix an element $y_{i}$ of $X$ on which the
displacement function of $f_{i}(g)$ attains its minimum
$\t(f_{i}(g))$, i.e.
$\t(f_{i}(g))=d\big(f_{i}(g)y_{i},y_{i}\big)=\textrm{inf\,}\{d\big(f_{i}(g)y,y\big)|\,y\in
X\}$.\footnote{The reader should not confuse this with the algebraic translation length of 
the elements of a group $G$.} 
Since $X$ is a $\delta$-hyperbolic space, there is a non-negative
constant $K(\delta)$ depending only on $\delta$ such that
\begin{equation}\label{ineq:a}
d\big(f_{i}(g)x_{i},x_{i}\big)\geq
2d(x_{i},y_{i})+\t(f_{i}(g))-K(\delta).
\end{equation}
Then,
\begin{equation} \label{ineq:b}
\begin{array}{ccl} A & = &
2d\big(f_{i}(g)x_{i},x_{i}\big)-d\big(f_{i}(g)^{2}x_{i},x_{i}\big)\\
 & \geq &
 4d(x_{i},y_{i})+2\t(f_{i}(g))-2K(\delta)-d\big(f_{i}(g)^{2}x_{i},x_{i}\big)\\
 & \geq & 4d(x_{i},y_{i})+2\t(f_{i}(g))-2K(\delta)-2d(x_{i},y_{i})-2\t(f_{i}(g))\\
 & = & 2d(x_{i},y_{i})-2K(\delta)\,,
 \end{array}\end{equation}
where the second inequality follows from the triangle one. Working
in a similar way, we see that
\begin{equation}
2d(x_{i},y_{i})-K(\delta)\leq
d\big(f_{i}(g)x_{i},x_{i}\big)-\t(f_{i}(g))\leq 2d(x_{i},y_{i})
\end{equation}
and hence we can say that 
\begin{equation}\label{ineq:d}
\big|d\big(f_{i}(g)x_{i},x_{i}\big)-\t(f_{i}(g))\big|\leq
2d(x_{i},y_{i})+K(\delta).
\end{equation}
Consequently,
\[\begin{array}{ccl}
\Big|\t_{\omega}(g)-\frac{\t(f_{i}(g))}{\lambda_{i}}\Big| & \leq &
\Big|\t_{\omega}(g)-d_{i}\big(f_{i}(g)x_{i},x_{i}\big)\Big|+
\Big|d_{i}\big(f_{i}(g)x_{i},x_{i}\big)-
\frac{\t(f_{i}(g))}{\lambda_{i}}\Big| \\
 & \leq &
 \Big|\t_{\omega}(g)-d_{i}\big(f_{i}(g)x_{i},x_{i}\big)\Big|+2\dsp\frac{d(x_{i},y_{i})}
{\lambda_{i}}+\frac{K(\delta)}{\lambda_{i}}\\
 & \leq &
 \Big|\t_{\omega}(g)-d_{i}\big(f_{i}(g)x_{i},x_{i}\big)\Big|+\dsp\frac{|A|}{\lambda_{i}}+
 3\frac{K(\delta)}{\lambda_{i}}\,,\\
\end{array}\]
where the second inequality follows from ~(\ref{ineq:d}) and the third
one from ~(\ref{ineq:b}). The $\omega$-limit of each term in
the right-hand side of the above inequality is $0$ (for the second term
see ~(\ref{eq:2})). Therefore
\[\t_{\omega}(g)=\textrm{lim}_{\omega}\frac{\t(f_{i}(g))}{\lambda_{i}}=0\,,\]
since $\t(f_{i}(g))=\t(g)$ for all $i$ ($f_{i}$ being a conjugating
automorphism). Hence, each element $g$ of the finitely generated
group $G$ acts as an elliptic isometry on $X_{\omega}$. This means
that the action of $G$ on $X_{\omega}$ has a global fixed point, say $z=(z_i)$. 
It follows that for every $\varepsilon>0$ and every
finite subset $F$ of $G$, there is an $\Omega\in \omega$ such that
\[d_{i}\big(z_{i},f_{i}(g)z_{i}\big)=\frac{d(z_{i},f_{i}(g)z_{i})}{\lambda_{i}}<
\varepsilon,\] for all $i\in \Omega$ and $g\in F$. This
contradicts to the minimality of $\lambda_{i}$ (see ~(\ref{eq:min})).
\end{proof}

We should mention here that this is the best possible result in that generality.
Indeed, as shown by Burnside \cite{burnside} and subsequently by several other
authors (see also Sah \cite{sah}), there are finite groups that posses non-trivial 
outer conjugating automorphisms
(known also as outer, class preserving automorphisms). Therefore one can 
easily construct free-by-finite groups with outer conjugating automorphisms by
considering the direct product of the above mentioned finite groups by free groups.  

We are now able to show our main theorem.

\begin{thm}\label{theorem}
Let $G$ be a conjugacy separable, hyperbolic group. Then each
virtually torsion-free subgroup of the outer automorphism group
$\out(G)$ of $G$ is residually finite.
\end{thm}
\begin{proof} It suffices to show that each torsion-free subgroup $H$
of $\out(G)$ is residually finite. We consider the following short
exact sequence
\[1 \lar \conj(G)\big/\inn(G)\stackrel{i}{\lar} \aut(G)\big/\inn(G)\stackrel{\pi}{\lar} 
\aut(G)\big/\conj(G) \lar 1\,.\]
 By Lemma ~\ref{lem:finite} the first term is finite. This implies that the 
restriction of $\pi$ on $H$
is a monomorphism. It follows that $H$ is residually finite being
isomorphic to a subgroup of $\aut(G)\big/\conj(G)$, which is
residually finite by Lemma ~\ref{lem:Gr}.
\end{proof}

In view of the above theorem, given a group $G$ it is natural to
seek conditions under which the outer automorphism group $\out(G)$
of $G$ is virtually torsion-free. Recently, Guirardel and Levitt ~\cite{gl} have shown that the
outer automorphism group of a hyperbolic group $G$ is virtually
torsion-free, provided that $G$ is virtually torsion-free. Thus, by
Theorem ~\ref{theorem}, the outer automorphism group of a
conjugacy separable hyperbolic group is residually finite, since 
a residually finite hyperbolic group is virtually torsion-free. 

The next lemma is more or less
known (see ~\cite{gl,Mc}).

\begin{lem}\label{lem:vir} Let $G$ be a finitely generated group containing a
normal subgroup $N$ of finite index whose center is trivial. If
$\out(N)$ is virtually torsion-free, then so is $\out(G)$.
\end{lem}

\begin{proof} Let $\aut_N(G)$ denote the subgroup of $\aut(G)$
consisting of those automorphisms of $G$ which fix $N$ and induce
the identity on $G/N$. The restriction map $\phi: \aut_N(G) \lar \aut(N)$ is an
injection. Indeed, let $g\in G$ and $f\in \aut_N(G)$. Then
$f(g)=gh_{g}$ for some $h_{g}\in N$. Suppose now that $f$ is in
the kernel of the restriction map. Then for each $h$ in $N$ we
have $ghg^{-1}=f(ghg^{-1})=gh_{g}hh_{g}^{-1}g^{-1}$. This implies
that $h_{g}$ is in the center of $N$, which is trivial, and
therefore $f$ is the identity. From the injectivity of $\phi$ we
get $\phi\big(\inn_N(G)\big)=\inn(N)$, where $\inn_N(G)$ denotes
the (normal) subgroup of $\inn(G)$ consisting of all inner
automorphisms of $G$ induced by elements of $N$. We conclude that
the quotient group $\aut_N(G)/\inn_N(G)$ embeds into $\out(N)$. In
particular, $\aut_N(G)/\inn_N(G)$ is virtually torsion-free.

Now the restriction of the natural projection $\pi:\aut(G) \lar
\out(G)$ to $\aut_N(G)$ induces a map $\tilde{\pi}:\aut_N(G)/\inn_N(G) 
\lar \out(G)$. The kernel of $\tilde{\pi}$ (which is equal to
$\inn(G)/\inn_N(G)$) is finite, since $N$ is of finite index in
$G$, while its image has finite index in $\out(G)$, since $\aut_N(G)$ 
is of finite index in $\aut(G)$, by ~\cite[Lemma 1]{Mc}. It
follows that each finite-index, torsion-free subgroup of $\aut_N(G)/\inn_N(G)$ 
embeds as a subgroup of finite index in $\out(G)$,
which proves the lemma.
\end{proof}

\begin{cor}\label{cor:free-by-finite}
The outer automorphism group of a finitely generated,
free-by-finite group is residually finite.
\end{cor}
\begin{proof} It is well-known that every finitely generated,
free-by-finite group $G$ is hyperbolic. Furthermore, $G$ is
conjugacy separable by ~\cite{Dy}. Therefore Theorem
~\ref{theorem} applies to $G$. On the other hand, it is also known
that the outer automorphism group of a free group is virtually
torsion-free. If $G$ is virtually infinite cyclic, then $\out(G)$
is finite. In the case where $G$ is not virtually cyclic, the
center of a free subgroup of finite index in $G$ is trivial, and
the result follows from Lemma ~\ref{lem:vir}.
\end{proof}

Recall that a Fuchsian group is a discrete subgroup of the group of 
isometries of the hyperbolic plane $\mathbb{H}^{2}$.

\begin{cor} The outer automorphism group of a finitely generated,
Fuchsian group is residually finite.
\end{cor}

\begin{proof} Every finitely generated Fuchsian group $G$ contains  
a normal torsion-free subgroup of finite 
index, say $N$, which is either a free group or the fundamental group of an
orientable surface group of genus $g\ge 2$. So $N$ is hyperbolic since it is either free or 
quasi-isometric to $\H^2$. Moreover $N$ is conjugacy separable \cite[Theorem 3.3]{stebe} and 
its outer automorphism group $\out(N)$ is virtually torsion-free and so $N$ 
satisfies the hypotheses of
Theorem \ref{theorem}. Hence, $\out(N)$ is residually finite.

If $N$ is cyclic then $G$ is virtually cyclic and so $\out(G)$ is finite. In all other
cases, $N$ is a normal subgroup of finite index in $G$ with trivial centre (since it is
a non-cyclic torsion-free hyperbolic group) and so from the proof of 
Lemma \ref{lem:vir} we have that $\aut_N(G)/\inn_N(G)$ is a subgroup of $\out(N)$. 
Hence, every subgroup of $\aut_N(G)/\inn_N(G)$ is residually finite. But, again from the proof of 
Lemma \ref{lem:vir}, there is a torsion-free subgroup of $\aut_N(G)/\inn_N(G)$ 
that  embeds as a finite index subgroup 
in $\out(G)$. Therefore, $\out(G)$ is residually finite.
\end{proof}

The corollary below generalizes the results of Grossman \cite{Gr} and of 
Allenby, Kim and Tang \cite{akt}. Notice that for the exceptional cases of the Torus and
the Klein bottle it is easily verified that the result still holds.
\begin{cor}
The mapping class group $\mod_S$ of a hyperbolic 2-orbifold $S$ with finite volume
is residually finite.
\end{cor}

\begin{proof}
The fundamental group of a hyperbolic 2-orbifold with finite volume is a
finitely generated Fuchsian group. Hence, the proof is an immediate consequence 
of the results of Fujiwara \cite{fujiwara}, the
above corollaries and the fact that subgroups of residually finite groups
are residually finite. 
\end{proof}

\subsubsection*{Addendum. Added, October 8, 2005.} Only recently the authors found out that their 
main result can be
generalized to relatively hyperbolic groups by using \cite[Theorem
1.1]{bs}.

Relatively hyperbolic groups were introduced by
Gromov in \cite{gromov}, in order to generalize notions such as
the fundamental group of a complete, non-compact, finite volume 
hyperbolic manifold and to give a hyperbolic version of small 
cancellation theory over free groups by adopting the geometric language 
of manifolds with cusps. 

This notion has been 
developed by several authors and, in particular, various 
characterizations of relatively hyperbolic groups
have been given (see \cite{bow,osin} and \cite{ds} 
and references therein).

We recall here one of Bowditch's equivalent definitions. A finitely
generated group $G$ is {\it hyperbolic relative to a
family of finitely generated subgroups\/} $\mathcal{G}$ if $G$ admits a 
proper, discontinuous and isometric action on a proper, hyperbolic 
path metric space $X$ such that $G$ acts on the ideal boundary of 
$X$ as a geometrically
finite convergence group and the elements of $\mathcal G$
are the  maximal parabolic subgroups of $G$.

 We should mention here that Farb
\cite{farb} introduced a weaker notion of relative hyperbolicity 
for groups using constructions on the Cayley graph of the groups.





Except of the fundamental groups of hyperbolic manifolds
of finite volume, another interesting family of relatively 
hyperbolic groups are the fundamental groups of graphs of finitely generated
groups with finite edge groups which are hyperbolic relative to the 
family of vertex groups (since their action on the Bass-Serre tree 
satisfies definition 2 in \cite{bow}).

The key lemma of the paper, Lemma 2.2, generalizes to
relatively hyperbolic groups.

\setcounter{section}{2} \setcounter{thm}{1}

\begin{lem}\hskip -.2cm{$\mathbf '$}
Let $G$ be a relatively hyperbolic group. Then the group
Inn$(G)$ of inner automorphisms of $G$ is of finite
index in Conj$(G)$.
\end{lem}

\noindent {\it Sketch of proof.} In this case the Cayley graph
of $G$ is replaced by the $\d$-hyperbolic metric space $X$ on
which $G$ acts by isometries.

Let $\l_i=\inf\limits_{x\in X}\max\limits_{s\in S} d(x,f_i(s)x)$
where $S$ is a fixed finite generating set $G$ and let $x_i^0\in X$
such that $\max\limits_{s\in S}d(x_i^0,f_i(s)x_i^0)\le \l_i+\frac 1i$.

By the proof of Theorem 1.1 in \cite{bs}, the sequence $\l_i$
converges to infinity. Hence for every non-primitive ultrafilter
$\w$ on $\N$ the based ultralimit $(X_{\w},d_{\w},x_{w}^0)$ of the
sequence of based metric spaces $(X,\frac{d}{\l_i},x_i^0)$ is an
$\R$-tree. Moreover, there is an induced non-trivial isometric
$G$-action on $(X_{\w},d_{\w},x_{w}^0)$, given by $g\cdot
x_i=f_i(g)x_i$.

Following step-by-step the proof of Lemma 2.2 we arrive again at the
contradiction that the action has a global fixed point.

The reader should be careful in the following.
\begin{enumerate}
\item The elements $y_i$ of $X$ are chosen such that
$$\t(f_i(g))\le d(f_i(g)y_i,y_i) \le \t(f_i(g))+\frac 1i.$$
\item The inequalities (5), (6) and (7) become
$$ A\ge 2d(x_i,y_i)-2K(\d) -\frac 2i \eqno (5') $$
$$ 2d(x_i,y_i)-K(\d) \le d(f_i(g)x_i,x_i)-\t(f_i(g))\le 2d(x_i,y_i)+\frac 1i \eqno (6')$$
$$ \textrm{and} \quad |d(f_i(g)x_i,x_i)-\t(f_i(g))| \le 2d(x_i,y_i)+K(\d)+\frac 1i, \eqno (7')$$
\noindent respectively.

\item Finally, the last inequality becomes
$$\left|\t_{\w}(g) - \frac{\t(f_i(g))}{\l_i}\right|\le
\left| \t_{\w}(g)- d_i(f_i(g)x_i,x_i)\right|+\frac{|A|}{\l_i}+3\frac{K(\d)}{\l_i}+\frac{3}{i\l_i}.$$
\end{enumerate}
\hfill$\Box$

Consequently, Theorem 2.3 generalizes as follows.

\setcounter{thm}{2}
\begin{thm}\hskip -.2cm{$\mathbf '$}
Let $G$ be a conjugacy separable, relatively hyperbolic group.
Then each virtually torsion-free subgroup of the outer
automorphism group {\rm Out}$(G)$ of $G$ is residually finite.
\end{thm}

\noindent Department of Mathematics, University of the Aegean, Karlovassi, 832 00 Samos, 
Greece {\it email: vmet@aegean.gr} 

\bigskip

\noindent Amalthias 18, 412 22 Larissa Greece {\it email: msikiot@cc.uoa.gr}\\

{\it Current address:} Department  of Mathematics and Statistics, University of Cyprus, 
P.O. Box 20537, 1678 Nicosia, Cyprus {\it email: msikiot@ucy.ac.cy} 
\end{document}